\documentstyle{amsart}
\title{$A$ -- sets.}
\author{Tomek Bartoszynski and Marion Scheepers}
\address{Department of Mathematics\\ Boise State University\\ Boise,
   Idaho 83725}
\thanks{The second author was supported by Idaho State Board of
   Education grant 94-051.}
\email{tomek@@math.idbsu.edu\\ marion@@math.idbsu.edu}
\newcommand{\reals}{\Bbb R}
\newcommand{\naturals}{\Bbb N}
\newcommand{\rationals}{\Bbb Q}

\newtheorem{th}{Theorem}
\newtheorem{cor}[th]{Corollary}

\newtheorem{proposition}[th]{Proposition}

\begin{document}
\maketitle
\vspace{0.5in}

   According to E. Borel \cite{B} a set $X$ of real numbers has {\em strong
     measure zero} if there is for each sequence
   $\epsilon_1,\epsilon_2,\dots,\epsilon_n,\dots$ of positive real
   numbers a corresponding sequence $I_1, I_2, \dots, I_n, \dots$ of
   open intervals such that $length(I_n)<\epsilon_n$ for each $n$, and
   $X\subseteq\cup_{n=1}^{\infty}I_n$. The
   collection of strong measure zero sets forms a $\sigma$--ideal.
   Galvin, Mycielski and Solovay \cite{G-M-S} discovered that
   the notion of a strong measure zero set has an algebraic
   characterization: $X$ has strong measure zero if, and only if,
   there is for each first category set $M$ a real number $t$ such
   that $X\cap(t+M)=\emptyset$.

   By analogy a set $Y$ was defined to have strong first category if
   there is for each measure zero set $N$ a real number $t$ such that
   $Y\cap(t+N)=\emptyset$. Some of the basic
   combinatorial properties of the collection of strong first category
   sets are not easily discernible from this algebraic definition:
   Only recently Pawlikowski proved that Sierpinski sets have
   strong first category. At present it is not even known if this
   collection of sets is an ideal!

   In this connection, the first author and H. Judah, pursuing ideas
   of I. Rec\mbox{\l}aw, defined and studied properties of the so--called
   $R^{\cal M} - $, $SR^{\cal M} - $, $R^{\cal N} - $ and $SR^{\cal
     N} - $ sets \cite{B-J1}. For convenience of the reader we briefly define 
   these notions here. ${\cal M}$ denotes the ideal of
   first category subsets of the real line, while the symbol ${\cal
     N}$ denotes the ideal of measure zero sets. 
    Suppose that ${\cal J}$ is an ideal of subsets of the real line.
    A Borel set $H \subseteq \reals \times \reals$ is called 
    ${\cal J}$-set  if $(H)_x =\{y:(x,y)\in H\}\in {\cal J} $ for all
     $ x \in \reals$.  We say that $X \subseteq \reals$ is an
     $R^{{\cal J}}$ set if  for every ${\cal J}$-set $H$,
   $\bigcup_{x \in X} (H)_x \neq \reals .$
   A set $X \subseteq \reals$ is an $SR^{{\cal J}}$ set if 
   for every ${\cal J}$-set $H$, $\bigcup_{x \in X} (H)_x \in {\cal J} .$

   Independently, the second author started
   searching for covering properties of the collection of strong first category
   sets, analogous to the properties studied by Menger \cite{M},
   Hurewicz \cite{Hu} and Rothberger \cite{Ro}.

   A common thread seems to be emerging from these two approaches. In
   this paper we report some of our findings to illustrate this. 
   We introduce what we call $A$ -- sets. The motivation for the name
   is simply that letters from the alphabet seem to be used for
   various small sets studied by various authors (such as $C$--sets,
   $Q$--sets, and so on) and the letter $A$ had not been reserved yet.

   For a set $X$ in the range of a function $f$, we write
   $f^{\leftarrow}[X]$ to denote the pre--image of $X$ under $f$.
   Where not specified, our notation or terminology follows that of 
   \cite{B-J1} or \cite{F-M}.

\section{The $A_1$-sets.}

   Let $X$ be a set of real numbers.
    $X$ is an $A_1$-set if: For every sequence $({\frak
   U}_1,{\frak U}_2,\dots)$ such that $X\subseteq\bigcup{\frak U}_n$ for
   each $n$, and each ${\frak U}_n$ is a {\em countable} collection of
   Borel sets, there is a sequence $(Y_1,Y_2,\dots)$
   such that for each $n$ we have $Y_n\in{\frak U}_n$ and
   $X\subseteq\cup_{n<\infty}Y_n$.

   The $A_1$--property is a direct analogue of Rothberger's
   property $C''$. Indeed, we have:

\begin{th}\label{rothbmeager} For a set $X$, the 
  following are equivalent: 
  \begin{enumerate}
    \item{$X$ is an $A_1$--set.}
    \item{Every Borel image of $X$ into $^{\naturals}\naturals$ has
        Rothberger's property $C''$.}
    \item{$X$ is an $R^{\cal M}$--set.}
  \end{enumerate}
\end{th} 

\begin{pf}
  $(1)\Rightarrow (2):$\ \ Let
  $\Psi:X\rightarrow\mbox{$^{\naturals}\naturals$}$ be a Borel
  function. Let $({\frak U}_1,{\frak U}_2,\dots,{\frak U}_n,\dots)$ be a
  sequence of open covers of $\Psi[X]$; we may assume that each is
  countable. Enumerate each ${\frak U}_n$ bijectively as $\{U^n_1,
  U^n_2,\dots, U^n_k,\dots\}$. 

  Then $\Psi^{\leftarrow}[U^n_k]=F^n_k$ is a Borel subset of $X$ for
  each $(n,k)$. Choose Borel subsets $A^n_k$ of $\reals$ such that
  $F^n_k=X\cap A^n_k$.
  For each $n$ put ${\cal S}_n=\{A^n_k:k\in{\naturals}\}$.
  Then each ${\cal S}_n$ is a cover of $X$ by countably many 
  Borel sets. Select for each $n$ a $k_n$ such that
  $X\subseteq\cup_{n=1}^{\infty}A^n_{k_n}$. Then
  $\Psi[X]\subseteq\cup_{n=1}^{\infty}U^n_{k_n}$. \\
  $(2)\Rightarrow(1):$\ \ Assume that every Borel image of $X$ into
  $^{\naturals}\naturals$ has property $C''$. Let $({\frak U}_1, {\frak
    U}_2,\dots, {\frak U}_n,\dots)$ be a sequence such that each ${\frak
    U}_n$ is a covering of $X$ by countably many Borel sets.
  Enumerate ${\frak U}_n$ bijectively as $\{U^n_1, U^n_2,\dots,
  U^n_k,\dots\}$. Define
\[\Psi:X\rightarrow\mbox{$^{\naturals}\naturals$}\]
  so that for each $x\in X$ we have $\Psi(x)(n)=\min\{m:x\in U^n_m\}$.
  Note that for every basic open set $[(n_1,\dots,n_k)]$,
  $\Psi^{\leftarrow}[[(n_1,\dots,n_k)]] = X\cap(U^1_{n_1}\cap\dots\cap
  U^k_{n_k})\setminus(U^1_{1}\cup\dots\cup U^1_{n_1-1} \cup\dots\cup
  U^k_1\cup\dots \cup U^k_{n_k-1})$, a Borel subset of $X$. Then
  $\Psi$ is a Borel mapping. By Hypothesis $\Psi[X]$ has property
  $C''$. Thus, we find a $g$ such that for each $x$ the set
  $\{n:\Psi(x)(n)=g(n)\}$ is infinite (see \cite{Re}, Lemma 1). But
  then $X\subseteq \cup_{n=1}^{\infty} U^n_{g(n)}$.

  As to the equivalence of $(2)$ and $(3)$:  According to Theorem 4.3
  of \cite{B-J1}, $X$ is an $R^{\cal M}$ -- set if, and only if, every
  Borel image of it into $^{\naturals}\naturals$ has the Rothberger property.
\end{pf}

   In particular, we see that every set which has property $A_1$ also
   has Rothberger's property $C''$. Since for example the Sierpinski
   sets do not have measure zero, while $C''$ -- sets do, we see that
   the collection of $A_1$ -- sets is too small to capture all the
   sets of strong first category.

\begin{cor}\label{a1.2} The collection of $A_1$--sets is a
  $\sigma$-ideal.
\end{cor}

\begin{pf}
  It is clear from the definition of these sets that this collection
  is closed under countable unions. Being an
  $R^{\cal M}$--set is hereditary, as is clear from the definition of
  that notion.
\end{pf}

   Since we have the inclusions $SR^{\cal N}\subseteq SR^{\cal
   M}\subseteq R^{\cal M}$, it follows that $SR^{\cal N}$-sets and $SR^{\cal
   M}$-sets are $A_1$-sets.

\section{The $A_2$--sets.}

   In his 1927 paper \cite{Hu} Hurewicz introduced a property of sets
   which is nowadays referred to as property $H$; see \cite{F-M}.
  A set $X$ of real numbers has property $H$ if for every sequence
  $\langle{\frak G}_n:n<\omega\rangle$ of open covers of $X$, there
  exists a sequence $\langle{\frak F}_n:n<\omega\rangle$ such that
  each ${\frak F}_n$ is a finite subset of ${\frak G}_n$, and such
  that $X\subseteq\cup_{n=1}^{\infty}(\cap_{m=n}^{\infty}\cup{\frak F}_m)$.

   Let $X$ be a set of real numbers. We say that $X$ is
   an $A_2$--set if:  For every sequence $({\frak
   U}_1,{\frak U}_2,\dots)$ such that $X\subseteq\cup{\frak U}_n$ for
   each $n$, and each ${\frak U}_n$ is a {\em countable} collection of
   Borel sets, there is a sequence $({\frak G}_1,{\frak G}_2,\dots)$
   such that for each $n$ we have ${\frak G}_n$ is a finite subset of
   ${\frak U}_n$ and 
   $X\subseteq\cup_{n=1}^{\infty}(\cap_{m=n}^{\infty}\cup{\frak G}_m)$.

   According to \cite{B-J1}, a set has property ${\cal H}$ if every
   Borel image of it in $^{\naturals}\naturals$ is bounded. As noted
   in \cite{B-J1}, every Sierpinski set is an $R^{\cal N}$ set with
   property ${\cal H}$. It follows that $R^{\cal N}\cap{\cal
     H}\not\subset A_1$.

   \begin{th}\label{hurewmeager}
     For a set $X$, the following are equivalent:
     \begin{enumerate}
       \item{$X$ is an $A_2$ - set.}
       \item{Every Borel image of $X$ into $^{\naturals}\naturals$ is
           bounded.}
       \item{Every Borel image of $X$ into $^{\naturals}\naturals$ has
           property $H$.}
     \end{enumerate}
   \end{th}

   \begin{pf}
     That $(1)$ implies $(2)$: Let $\Psi:X\rightarrow$
     $^{\naturals}\naturals$ be a Borel function. Let ${\cal O}_n =
     \{[\sigma]: \sigma\in\mbox{$^{n}\naturals$}\}$ for each $n$. Then $({\cal
       O}_1,{\cal O}_2,\dots,{\cal O}_n,\dots)$ is a sequence of
     open covers of $\Psi[X]$. Enumerate
     each ${\cal O}_n$ bijectively as
     $\{[\sigma^n_1],[\sigma^n_2],\dots,[\sigma^n_k],\dots\}$. Then
     $\Psi^{\leftarrow}[[\sigma^n_k]] (\stackrel{def}{=} F^n_k$, say)
     is a Borel subset of 
     $X$ for each $(n,k)$. For each $(n,k)$ choose a Borel subset
     $B^n_k$ of $\reals$ such that $F^n_k=X\cap B^n_k$.

     For each $n$, put ${\frak U}_n=\{B^n_k:k\in\naturals\}$.
     Then $({\frak U}_1,{\frak U}_2,\dots,{\frak U}_n,\dots)$ is a
     sequence of countable Borel covers of $X$. Since
     $X$ is an $A_2$-set, we find for each $n$ a finite set ${\cal
       V}_n\subset{\frak U}_n$ such that
     $X\subseteq\cup_{n=1}^{\infty}(\cap_{n=m}^{\infty}\cup{\cal
       V}_n)$. We may assume that there is for each $n$ a $k_n$ such
     that ${\cal V}_n=\{B^n_j:j\leq k_n\}$. 

     Putting $g(n)=\max\{\sigma^{n+1}_j(n):j\leq k_{n+1}\}$, we see
     that $\Psi(x)\prec g$ for each $x\in X$.\\

   To see that $(2)$ implies $(3)$, consider any Borel function $\Psi$
   from $X$ to $^{\naturals}\naturals$. For any continuous function
   $h$ from $\Psi[X]$ to $^{\naturals}\naturals$, the function
   $h\circ\Psi$ is a Borel function on $X$, so that $h\circ\Psi[X]$ is
   bounded. But this shows that every continuous image of $\Psi[X]$ is
   bounded; by \S 5 of \cite{Hu}, $\Psi[X]$ has property $H$.\\
 
     To see that $(3)$ implies $(1)$, let $({\frak U}_1, {\frak
       U}_2,\dots,{\frak U}_n,\dots)$ be a sequence such that for each
     $n$, ${\frak U}_n$ is a countable cover of $X$ with Borel
     sets. Enumerate each ${\frak U}_n$ bijectively as $\{U^n_1,
     U^n_2,\dots, U^n_k,\dots\}$. For each $x\in X$, define
\[\Psi(x)(m)=\min\{n:x\in U^m_n\}.\]
     Then $\Psi$ is a Borel mapping, and so $\Psi[X]$ has property
     $H$. Since the identity mapping is continuous, $\Psi[X]$ is a
     bounded subset of $^{\naturals}\naturals$. Choose $g$ such that
     $\Psi(x)\prec g$ for each $x$. Put ${\cal V}_n=\{U^n_i:i\leq
     g(n)\}$ for each $n$. Then
     $X\subseteq\cup_{n=1}^{\infty}(\cap_{m=n}^{\infty}\cup{\cal
       V}_m)$.
\end{pf}
   
   We see that $R^{\cal N}\cap{\cal H}\subset A_2$; in particular,
   every Sierpinski set is an $A_2$--set.

\section{The $A_3$--sets.}

   In his 1924 paper \cite{M}, Menger introduced a property of sets
   which is related to the Rothberger property. This Menger property
   is nowadays referred to as property $M$ (see \cite{F-M}). 
   A set $X$ of real numbers has property $M$ if for every sequence
   $\langle{\frak G}_n:n<\omega\rangle$ such that each ${\frak G}_n$
   is a chain of open sets covering $X$, there
   exists a sequence $\langle F_n:n<\omega\rangle$ such that
   each $F_n\in{\frak G}_n$, and such that $\{F_n:n<\omega\}$ is an
   open cover of $X$. 

   The proof of Theorem 3 on p. 21 of \cite{F-M} shows, {\em mutatis
   mutandis} (see also \cite{Re}, Proposition 3), that

\begin{th}
  For a subset $X$ of the real line, the following are equivalent:
  \begin{enumerate}
    \item{$X$ has property $M$.}
    \item{No continuous image of $X$ into $^{\omega}\omega$ is a
        dominating family.}
  \end{enumerate}
\end{th}

   We now introduce our third covering property. Let $X$ be a 
   first category set of real numbers. 
    $X$ is an $A_3$-set if: For every sequence $({\frak
   U}_1,{\frak U}_2,\dots)$ such that $X\subseteq\cup{\frak U}_n$ for
   each $n$, and each ${\frak U}_n$ is a countable chain of 
   Borel sets, there is a sequence $(Y_1,Y_2,\dots)$
   such that for each $n$ we have $Y_n\in{\frak U}_n$ and
   $X\subseteq\cup_{n<\infty}Y_n$.
   We denote the collection of $A_3$--sets by the symbol ${\cal A}_3$.

  The collection ${\cal A}_3$ is closed under countable unions. 
  
\begin{th}\label{mengermeagerth}
     For $X$ a set of real numbers, consider the following
     assertions: $(1)$, $(2)$, $(3)$ and $(4)$ are
     equivalent and each implies $(5)$.
  \begin{enumerate}
    \item{$X$ has property $A_3$.}
    \item{For every sequence $({\frak U}_1, {\frak U}_2, \dots, {\frak
          U}_n, \dots)$ such that each ${\frak U}_n$ is a collection of
        countably many Borel sets whose union contains
        $X$, there is a sequence $({\cal V}_1, {\cal V}_2, \dots,
        {\cal V}_n,\dots)$ such that
        \begin{itemize}
          \item{${\cal V}_n$ is a finite subset of ${\frak U}_n$ for
              each $n$, and}
          \item{$X\subseteq\cup_{n=1}^{\infty}(\cup{\cal V}_n)$.}
        \end{itemize}}
    \item{Every Borel image of $X$ in $^{\omega}\omega$ has
        property $M$.}
    \item{No Borel image of $X$ in $^{\omega}\omega$ is dominating.}
    \item{$X$ has property $M$.}
   \end{enumerate}
\end{th}

   \begin{pf} The proof that $(1)$ and $(2)$ are equivalent is standard.
     To see that $(1)$ implies $(3)$, let $X$ be a set with property $A_3$,
     and let $\Psi$ be a Borel function from $X$ to
     $^{\omega}\omega$. We must show that $\Psi[X]$ has property $M$.

    Take a sequence $({\frak G}_n:n<\omega)$ of open covers of $\Psi[X]$. We
    may assume that each ${\frak G}_n$ is a countable ascending chain
    of open sets. Enumerate ${\frak G}_n$ in ascending order as
    $\{G^n_1,G^n_2,\dots,G^n_m,\dots\}$.

   For each $(n,m)$, put $F^n_m=\Psi^{\leftarrow}[G^n_m]$, a Borel
   subset of $X$. Observe that $F^n_m\subseteq F^n_{m+1}$
   for all $n$ and $m$. For each $(n,m)$ we fix a Borel subset
   $A^n_m$ of the real line such that $F^n_m=X\cap A^n_m$. 
   We may assume that $A^n_m\subseteq A^n_{m+1}$ for all $n$ and $m$.

   Define: ${\frak U}_n=\{A^n_1,A^n_2,\dots,A^n_k,\dots\}$. Then each
   ${\frak U}_n$ is an ascending chain of Borel sets,
   and $X\subseteq\cup{\frak U}_n$ for each $n$.

   Since $X$ is a $A_3$--set, we find for each $n$ a $k_n$ such that
   $X\subseteq\cup_{n=1}^{\infty}A^n_{k_n}$. But then
   $\Psi[X]\subseteq\cup_{n=1}^{\infty}G^n_{k_n}$, and we succeeded in
   finding the required selector for the given sequence of open
   covers.

   To see that $(3)$ implies $(1)$, let $({\frak U}_n:n<\omega)$ be a
   sequence of ascending chains of Borel sets, such
   that for each $n$ we have $X\subseteq\cup{\frak U}_n$. For each $n$
   enumerate ${\frak U}_n$ in ascending order as $\{U^n_1,
   U^n_2,\dots, U^n_k,\dots\}$. Now define 
\[\Psi:X\rightarrow\mbox{$^{\naturals}\naturals$}\]
   so that $\Psi(x)(n)=\min\{m:x\in U^n_m\}$. Since for a finite
   sequence $(i_1,\dots,i_n)$ of natural numbers we have
\[X\cap(U^1_{i_1}\cap\dots\cap U^n_{i_n})=\{x\in X:j\leq n\Rightarrow
\Psi(x)(j)=i_j\},\]
   we see that the inverse image of any basic open subset of
   $^{\naturals}\naturals$ is a Borel subset of $X$. Thus $\Psi$ is
   a Borel mapping. 

   Then $\Psi[X]$ has property $M$. As subset of
   $^{\naturals}\naturals$, it is not a dominating family. Choose a
   $g$ which is not dominated by the family $\Psi[X]$. For each $x\in
   X$ we see that there is an $n$ such that $\Psi(x)(n)\leq g(n)$.
   Finally we set $V_n=U^n_{g(n)}$ for each $n$. Then $V_n\in{\frak
     U}_n$, and $X\subseteq\cup_{n=1}^{\infty}V_n$.

   To see that $(3)$ implies $(4)$, note that if $\Psi$ is a Borel
   mapping and $g$ is continuous, then $g\circ\Psi$ is a Borel
   mapping, and the identity mapping $I(x)=x$ of
   $^{\omega}\omega$ is continuous. Thus, if a Borel image of
   $X$ is dominating, we would have a subset of $^{\omega}\omega$ which
   has property $M$ and yet is dominating, a contradiction.

   That $(4)$ implies $(3)$ can be seen as follows: consider a Borel
   image of $X$. It is not dominating. Following this Borel
   function with a continuous function results in a Borel
   function, and so this second image is still not
   dominating. Thus, no continuous image of the image of a Borel
   mapping of $X$ into $^{\omega}\omega$ is dominating. 

   To see that $(4)$ implies $(5)$, we simply note that continuous
   functions are Borel functions; thus no continuous image of
   $X$ is dominating. But this implies that $X$ has property $M$.
\end{pf}

\begin{cor}\label{a2.1} Every $A_1$--set and every $A_2$--set is an $A_3$ set.
\end{cor}

\begin{cor}\label{a2.2} Property $A_3$ is hereditary, ${\cal A}_3$ is
  a $\sigma$-ideal.
\end{cor}

\begin{pf}
  The property is hereditary on account of $(4)$ of Theorem
  \ref{mengermeagerth}. Since $^{\omega}\omega$ is a Borel image of
  $\reals$, it also follows from Theorem \ref{mengermeagerth} that
  ${\cal A}_3$ is a $\sigma$-ideal.
\end{pf}

   If $X\subset Z$ is a subset of a
   topological space $Z$, and if $D\subset Z$ is countable, then any
   continuous function $\pi:X\rightarrow Y$ can be extended to an
   $F_{\sigma}$ (and {\em ipso facto} a Borel) function $\rho:X\cup
   D\rightarrow Y$. 

   In \cite{W1} Rec\mbox{\l}aw shows that $MA$ implies the existence of a set
   of real numbers $X$ with the properties (see his Theorem 3):
   \begin{itemize}
     \item{There exists a continuous function from $X$ onto the closed
         unit interval,}
     \item{$X+F$ has Lebesgue measure zero whenever $F$ does, and}
     \item{there is a countable set $D$ such that $X\cup D$ is a
         $\gamma$--set.} 
   \end{itemize}

   This set $X$ cannot be an $A_3$--set. For the set
   $Y=[0,1]\setminus\rationals$ is homeomorphic with $^{\omega}\omega$
   (via continued fraction expansion), and this homeomorphism can be
   extended to all of $[0,1]$ in such a way that the resulting map is
   $F_{\sigma}$. Then we have an $F_{\sigma}$ function from $X$ onto
   $^{\omega}\omega$. Indeed, we find an $F_{\sigma}$ map from $X\cup
   D$ onto $^{\omega}\omega$. Accordingly, null--additive (and thus
   first category--additive) sets need not be $A_3$--sets. In particular,
   strong first category sets or even first category sets of strong
   measure zero need not be $A_3$ sets. Also, $\gamma$--sets
   need not be $A_3$--sets.

\begin{proposition}\label{Badd} $add({\cal A}_3)={\frak b}$.
\end{proposition}

\begin{pf} It is clear that $add({\cal A}_3)$ is regular and
   uncountable. Let $\lambda<{\frak b})$ be a regular infinite cardinal
   number, and assume that for each $\mu<\lambda$ it had already been
   established that $\mu<add({\cal A}_3)$. Let $X_{\alpha}, \
   \alpha<\lambda$ be a sequence of elements of ${\cal A}_3$ such that
   $X_{\alpha}\subset X_{\beta}$ whenever $\alpha<\beta$. Put
   $X=\cup_{\alpha<\lambda}X_{\alpha}$. Let
   $({\frak U}_1,{\frak U}_2,{\frak U}_3,\dots,{\frak U}_n,\dots)$ be a
   sequence such that each ${\frak U}_n$ is an ascending chain of Borel
   sets whose union contains $X$. For each $n$ we write
   ${\frak U}_n=\{U^n_1,U^n_2,\dots,U^n_m,\dots\}$ where $U^n_m\subset
   U^n_k$ whenever $m<k$.

   For each $\alpha<\lambda$, define a function $f_{\alpha}$ such that
   $X_{\alpha}\subseteq\cup_{n=1}^{\infty}U^n_{f_{\alpha}(n)}$. Since
   we have $\lambda<{\frak b}$, select a function $f$ such that
   $f_{\alpha}\prec f$ for each $\alpha<\lambda$.

   Next, fix $I\in[\lambda]^{\lambda}$ and $N<\infty$ such that 
\begin{enumerate}
\item{$m>N\Rightarrow f_{\alpha}(m)<f(m)$ for each $\alpha\in I$, and}
\item{$f_{\alpha}\lceil_{N+1}=f_{\beta}\lceil_{N+1} = \sigma$ for all
   $\alpha,\beta\in I$.}
\end{enumerate}

   Put $M_1=U^1_{\sigma(1)},\dots, M_N=U^N_{\sigma(N)}$, and
   $M_k=U^k_{f(k)}$ for each $k>N$. Then
   $X_{\alpha}\subseteq\cup_{k=1}^{\infty}M_k$ for each $\alpha\in I$,
   whence $X\subseteq\cup_{k=1}^{\infty}M_k$.
   This shows that ${\frak b}\leq add({\cal A}_3)$.

   To see that $add({\cal A}_3)\leq {\frak b}$ is easy, and left to
   the reader.
\end{pf}

   Thus, $MA$ implies that every set of reals of cardinality less than
   ${\frak c}$ is a $A_3$-set. There are first category sets which do
   not have property $A_3$.


\begin{thebibliography}{}

\bibitem{B-J1} T. Bartoszynski and H. Judah, {\em Borel images of sets
  of reals}, preprint.

\bibitem{B} E. Borel, {\em Sur la classification des ensembles de
   mesure nulle}, {\bf Bulletin de la Societe Mathematique de France}
   47 (1919), 97 -- 125. 

\bibitem{F-M} D.H. Fremlin and A.W. Miller, {\em On some
  properties of 
  Hurewicz, Menger, and Rothberger}, {\bf Fundamenta Mathematicae} 129
  (1988), 17 -- 33.

\bibitem{G-M-S} F. Galvin, J. Mycielski and R.M. Solovay, {\bf Notices
  of the American Mathematical Society} 26 (1979), A -- 280.

\bibitem{Hu} W. Hurewicz, {\em \"Uber Folgen stetiger Funktionen},
   {\bf Fundamenta Mathematicae} 9 (1927), 193--204.

\bibitem{M} M. K. Menger, {\em Einige \"Uberdeckungss\"atze der
  Punktmengenlehre}, {\bf Sitzungsberichte der Wiener Akademie}, 133
  (1924), 421--444.

\bibitem{Ro} F. Rothberger, {\em Eine Versch\"arfung der Eigenschaft
  C}, {\bf Fundamenta Mathematicae} 30 (1938), 50 - 55.

\bibitem{W1} I. Rec\mbox{\l}aw, {\em On small sets in the sense of measure
   and category}, {\bf Fundamenta Mathematicae} 133 (1989), 255--260.

\bibitem{Re} I. Rec\mbox{\l}aw, {\em Every Lusin set is undetermined
  in Point-open game}, {\bf Fundamenta Mathematicae}, to appear.

\end{thebibliography}
\end{document}